\documentclass[12pt,a4paper,]{article}
\usepackage[latin1]{inputenc}
\usepackage[T1]{fontenc}
\usepackage[english]{babel}
\usepackage{amsmath,amssymb,amsfonts,amsthm}
\usepackage[pdftex]{graphicx}
\usepackage[pdftex]{color}
\usepackage{enumerate}
\usepackage{dsfont}
\usepackage[all]{xy}
\usepackage{fancybox}
\usepackage{makeidx}
\usepackage{multicol}
\usepackage{mathrsfs}
\usepackage{titlesec}
\usepackage{hyperref}
\usepackage{anysize}
\marginsize{2.45cm}{2.35cm}{1.6cm}{2.3cm} 

\makeindex

\newtheorem{teo}{Theorem}

\newtheorem{conj}[teo]{Conjecture}

\theoremstyle{definition}
\newtheorem{deff}[teo]{Definition}
\newtheorem{exa}[teo]{Example}
\newtheorem{rmk}[teo]{Remark}

\newcommand{\bs}{ \backslash }
\newcommand{\ce}{ \mathcal{E}}
\newcommand{\dd}{{\ddabc}}
\DeclareMathOperator{\ddabc}{d}
 \newcommand{\ddd}{{\dddabc}}
\DeclareMathOperator{\dddabc}{D}

\newcommand{\Nset}{\mathds{N}}

\newcommand{\Zset}{\mathds{Z}}

\newcommand{\Rset}{\mathds{R}}

\title{A bounded coarse structure for families of pseudometrics}
\author{Jesús P. Moreno-Damas}
\date{}

\begin{document}

\maketitle

\begin{abstract}We define the bounded coarse structure attached to a family
of pseudometrics and give some counterexamples to conjectures that arise naturally.\end{abstract}

\section{Introduction}
We define the bounded coarse structure attached to a family of
pseudometrics of a space, generalizating the one attached to a
single metric. We do it in the same way as it is done in
\cite{moreno1} for the $C_0$ coarse structure (defined originally by
Wright for a single metric in \cite{wright,wright3}).

Moreover, we give some counterexamples of conjectures which arise naturally.

They are some of the conclusions of a coarse geometry's seminar with
Jerzy Dydak, Manuel Alonso Morón,  and Jesús P. Moreno-Damas in the
Dept. of Geometría y Topología of Universidad Complutense de Madrid.

\section{Preliminaries}

A pseudometric on a set is a map $\dd:X\times X\rightarrow [0,\infty)$ such that
$\dd(x,x)=$ for every $x$, $\dd(x,y)=\dd(y,x)$ for every $x,y$ and
$\dd(x,z)\leq \dd(x,y)+\dd(y,z)$ for every $x,y,z$. If moreover
$\dd(x,y)=0$ implies $x=y$ for every $x,y$, then $\dd$is a metric.

\section{The bounded coarse structure attached to a family of pseudometrics}

If $\{\ce_i\}_{i\in I}$ is a family of coarse structures over $X$,
then $\bigcap_{i\in I}\ce_i$ is a coarse structure.

If particular, if $\ddd$ is a family of pseudometrics of $X$ and
$\ce_b(\dd_i)$ is the bounded coarse structure of each one, then
\begin{equation}\label{defcebpseudometricas}
\bigcap_{\dd\in \ddd}\,\,\,\ce_b(\dd)=\{E\subset X\times X:\forall \dd\in
\ddd\exists R_\dd>0\textrm{ such that } \dd(x,y)\leq R_\dd\forall
(x,y)\in E\}
\end{equation}
is a coarse structure.

\begin{deff}\label{definition1} Let $\ddd$ a family of pseudometrics of $X$. The
bounded coarse structure attached to $\ddd$, denoted by
$\ce_b(\ddd)$ is the one defined in
(\ref{defcebpseudometricas}).\end{deff}

\begin{rmk} Let $X$ be a set, let $\ddd$ be a family of
pseudometrics of $X$. Consider the metric on $X$:
$$\rho(x,y)=\left\{\begin{array}{ll}0&\textrm{ if }x=y\\1&\textrm{
if }x\neq y\end{array}\right.$$ and let
$\ddd'=\{\dd+\rho:\dd\in\ddd\}$.

Then, $\ddd'$ is a family of metrics of $X$ such that
$\ce_b(\ddd)=\ce_b(\ddd')$.\end{rmk}

The first natural question is if the bounded coarse structure
attached to a family of pseudometrics is more general than the one
attached to one metric. Example \ref{counterexample1} will show
that, in general, it doesn't happen.

\begin{deff}Let $(X,\ce)$ be a coarse space and let
$\{E_\lambda\}\subset \ce$. We say that $\{E_\lambda\}$ generes
strongly $\ce$ if for every $E\in\ce$, there exists
$\lambda\in\Lambda$ such that $E\subset E_\lambda$.\end{deff}

\begin{exa}\label{ejmetric} Let $(X,\dd)$ be a pseudometric space and let, for every
$n\in\Nset$:
$$E_n=\{(x,y):\dd(x,y)\leq n\}$$
Then, $\{E_n\}_{n=1}^\infty$ generes strongly $(X,\dd)$.\end{exa}

\begin{exa}\label{counterexample1} For each $n\in\Nset$, consider in $\Rset^\Nset$ the
pseudometric $\dd_n(x,y)=|x_n-y_n|$. Let $\ddd=\{\dd_n:n\in\Nset\}$.
From the definition, we get that $E\in\ce_b(\ddd)$ if and only if
there are $\{R_n\}_{n\in\Nset}$ such that $|x_n-y_n|\leq R_n$ for
every $(x,y)\in E$ and every $n\in\Nset$.

There is not a pseudometric $\dd$ on $\Rset^\Nset$ such that
$\ce_b(\ddd)=\ce_b(\dd)$. Suppose, on the contrary, there is one.
Then, by Example \ref{ejmetric}, there exists
$\{E_n\}_{n=1}^\infty\subset \ce_b(\ddd)$ which generes strongly
$\ce_b(\ddd)$.

For each $n$, $E_n\in \ce_b(\ddd)$. Let $R_n$ be such that
$|x_i-y_i|<R_n$ for every $(x,y)\in E_n$.  Consider the set
$E=\big\{(x,y)\in \Rset^\Nset:|x_n-y_n|\leq 2R_n\forall
n\in\Nset\big\}$. Let $x=(0)_{n\in\Nset}$ and $y=(2
R_n)_{n\in\Nset}$. Clearly, $(x,y)\in E$, but, for every
$n\in\Nset$, $(x,y)\not\in E_n$. Then, $E\not\subset E_n$ for every
$n$ and $\{E_n\}_{n\in\Nset}$ for every $n\in\Nset$.\end{exa}

But in Example \ref{counterexample1}, in general the controlled sets
are not proper, as they use to be in the usual examples of bounded
coarse structures attached to metric spaces, like the one attached
to a finite generated group, for example $\Zset$ of $\Zset^2$. What
happen if the coarse space $\ce_b(\ddd)$ is proper? We put the
question in the following terms:

\begin{conj} Suppose that $\ddd$ is a family of pseudometrics of $X$
such that the topology induced by $\ddd$ is locally compact and the
coarse structure $\ce_b(\ddd)$ induced by $\ddd$ is proper. Then,
$(X,\ce_b(\ddd))$ is coarsely metrizable, i. e, there exists $\dd$
such that $\ce_b(\ddd)=\ce_b(\dd)$.\end{conj}

The following counterexample will show that this conjecture is not true:

\begin{exa} Consider on $\Nset\times\Nset$ the maps:
$$f_0:\Nset\times\Nset\rightarrow\Nset, f_0(a,b)= 2^a$$
$$\forall n\in\Nset, f_n:\Nset\times\Nset\rightarrow\Nset, f_n(a,b)=\left\{\begin{array}{ll}b &\textrm{ if }a=n\\0&\textrm{ if }a\neq n\end{array}\right.$$
and, for each $n\in\Nset\cup\{0\}$, let $\dd_n$ be the metric such
that $\dd_n(x,y)=|f_n(x)-f_n(y)|$.

Let $\ce=\ce_b(\{\dd_n\}_{n=0}^\infty)$. Let us see that $\ce$ is a
proper coarse structure such that there is no pseudometric $\dd$ of
$X$ with $\ce=\ce_b(\dd)$.

Let $M\in\Nset$ and let $\{R_n\}_{n=1}^\infty\subset\Rset^+$. Let
$E_{M,\{R_n\}_{n=1}^\infty}$ be the set of those
$((a,b),(a',b'))\in(\Nset\times\Nset)\times (\Nset\times\Nset)$
satisfying any of this properties: \begin{itemize}
\item $a,b,a',b'\leq M$
\item $a=a'$ and $|b-b'|\leq R_a$
\end{itemize}

Let us see first that
$$\{E_{M,\{R_n\}_{n=1}^\infty}:M\in\Nset,R_n\in\Rset^+\forall
n\in\Nset\}$$ generes strongly $\ce$.

Fix $M\in \Nset$ and $\{R_n\}_{n=1}^\infty\subset \Rset^+$. Let
$R'_0=2^{M+1}$ and $R'_n=R_n+2M$ for every $n\in\Nset$.

Let $((a,b),(a',b'))\in E_{M,\{R_n\}_{n=1}^\infty}$. Suppose that
$a=a'$ and $|b-b'|\leq R_n$. Clearly,
$$\dd_a((a,b),(a',b'))=|b-b'|\leq R_a\leq R'_a$$
$$\dd_0((a,b),(a',b'))=|2^a-2^a|=0<R'_0$$
$$\dd_n((a,b),(a',b'))=|0-0|=0<R_n\forall n\in\Nset\bs\{a\}$$

Suppose now that  $a,a',b,b'\leq M$. Then,
$$\dd_0((a,b),(a',b'))=|2^a-2^{a'}|\leq 2^a+2^{a'}\leq 2^M+2^M=R_0$$
$$\dd_n((a,b),(a',b'))=|0-0|=0 \,\,\forall n\in\Nset\bs\{a,a'\}$$
If $a=a'$, then:
$$\dd_a((a,b),(a',b'))=|b-b'|\leq b+b'=2M\leq R'_a$$
And, if $a\neq a'$, then:
$$\dd_a((a,b),(a',b'))=|b-0|= b\leq M\leq R'_a$$
$$\dd_{a'}((a,b),(a',b'))=|0-b'|= b'\leq M\leq R'_{a'}$$
Then, $E_{M,\{R_n\}_{n=1}^\infty}\in \ce$

Suppose now that $E\in \ce$. Let
$\{R_n\}_{n=0}^\infty\subset\Rset^+$ such that
$\dd_n((a,b),(a',b'))\leq R_n$ for every $((a,b),(a',b'))\in E$ and
every $n\in\Nset\cup\{0\}$. Let $M_0\in\Nset$ be such that $M_0\geq
R_0$ and let $M\in\Nset$ be such that $M\geq M_0$ and $M\geq R_n$
for every $n\leq M_0$. Let $((a,b),(a',b'))\in E$.

If $a=a'$, then $|b-b'|=|f_a(a,b)-f_a(a',b')|\leq R_a$ and
$((a,b),(a',b'))\in E_{M,\{R_n\}_{n=1}^\infty}$.

Suppose now that $a\neq a'$. If $a>a'$, then
$\dd_0((a,b),(a',b'))=|2^a-2^{a'}|\geq |2^a-2^{a-1}|=2^{a-1}\geq a>
a'$. By symmetry, if $a<a'$, then $\dd_0((a,b),(a',b'))\geq a'>a$.
Hence:
$$a,a'\leq\dd_0((a,b),(a',b'))\leq
R_0\leq M_0\leq M$$

Since $a,a'\leq M_0$, we have:
$$M\geq R_a\geq
\dd_a((a,b),(a',b'))=|f_a(a,b)-f_a(a',b')|=|b-0|=b$$
$$M\geq R_{a'}\geq
\dd_{a'}((a,b),(a',b'))=|f_{a'}(a,b)-f_{a'}(a',b')|=|0-b'|=b'$$

Then, $a,a',b,b'\leq M$ and $((a,b),(a',b'))\in
E_{M,\{R_n\}_{n=1}^\infty}$. Thus, $E\subset
E_{M,\{R_n\}_{n=1}^\infty}$.

Let us see that $\ce$ is a proper coarse structure. Obviously, the
diagonal is an open neighborhood of the diagonal. Let $E\in\ce$ be
symmetric. To see that $E$ is proper, it is enough to see that
$E_{(a',b')}$ is finite for every $(a',b')\in \Nset\times\Nset$.

Let $M\in\Nset$ and $\{R_n\}_{n=0}^\infty$ such that $E\subset
E_{M,\{R_n\}_{n=1}^\infty}$. Let $(a,b)\in E_{(a',b')}$. Then,
$((a,b),(a',b'))\in E_{M,\{R_n\}_{n=1}^\infty}$. If $a,a',b,b'\leq
M$, then $(a,b)\in\{(x,y):x,y\leq M\}$. If $a=a'$ and
$|b-b'|\leq R_{a'}$, then $(a,b)\in\{(x,y):x=a',b'-R_{a'}\leq x\leq
b'-R_{a'}\}$, Then,
$$E_{(a',b')}\subset \{(x,y):x,y\leq M\}\cup
\{(x,y):x=a',b'-R_{a'}\leq x\leq b'-R_{a'}\}$$ which is finite.
Then, $E$ is proper and $\ce$ is proper.

Suppose that $\dd$ is a pseudometric on $X$ such that
$\ce=\ce_b(\dd)$. Then, by Example \ref{ejmetric}, there is
$\{E_m\}_{m\in\Nset}\subset \ce$ generating strongly $\ce$.

For every $m$, let $M^m\in\Nset$ and $\{R^m_n\}_{n=1}^\infty$ be
such that $E_m\subset E_{M^m,\{R^m_n\}_{n=1}^\infty}$.

For every $n\in\Nset$, let $R_n=R^n_n+1$. Let
$E=E_{1,\{R_n\}_{n=1}^\infty}$.

For every $m\in\Nset$, consider the point
$x_m=((m,R^m_m+M_m+2),(m,M_m+1))$.

We have that $m=m$ and $|R^m_m+M_m+2-(M_m+1)|=R^m_m+1=R_m$, hence
$x_m\in E$.

In the other hand, $m,R^m_m+M_m+2,m,M_m+1\not\leq M_m$ and $m=m$ but
$|R^m_m+M_m+2-(M_m+1)|=R^m_m+1\not<R^m_m$, hence $x_m\not\in E_m$.

Then, $E\not\subset E_m$ for every $m$ and hence, there is no $\dd$
such that $\ce=\ce_b(\dd)$.

\end{exa}

\end{document}